\documentclass{amsart}
\usepackage{amsmath, amssymb}
\usepackage[all]{xypic}

\newtheorem{theorem}[equation]{Theorem}

\newtheorem{proposition}[equation]{Proposition}
\newtheorem{corollary}[equation]{Corollary}

\numberwithin{equation}{section}

\newenvironment{example}
	{\refstepcounter{equation}\medskip\noindent{\bf Example \theequation.}}
	{\medskip}

\newcommand{\Hom}{{\rm Hom}}
\renewcommand{\hom}{{\rm hom}}

\newcommand{\End}{{\rm End}}

\newcommand{\Aut}{{\rm Aut}}

\newcommand{\zed}{{\mathbb Z}}

\newcommand{\shift}{\mathcal{S}}

\newcommand{\ang}[1]{\langle #1 \rangle}

\newcommand{\blank}{\mbox{$\underline{\makebox[10pt]{}}$}}

\DeclareMathOperator{\rgr}{gr-\!}
\DeclareMathOperator{\lgr}{\!-gr}
\DeclareMathOperator{\rmod}{mod-\!}
\DeclareMathOperator{\lmod}{\!-mod}
\DeclareMathOperator{\rmodu}{mod^{u}-\!}
\DeclareMathOperator{\lmodu}{\!-mod^{u}}

\newcommand{\st}{\left\vert\right.}

\newcommand{\bbar}[1]{\overline{#1}}

\newcommand{\Hfunct}[3]{H_{#1}({#2},{#3})}

\newcommand{\Gaut}{\Aut_*}

\newcommand{\m}{\circ}
\newcommand{\T}{{\mathcal T}}

\title{$G$-algebras, twistings, and equivalences of graded categories}
\author{Susan  J. Sierra}
\address{Department of Mathematics \\
University of Michigan\\
Ann Arbor, MI 48109}
\email{ssierra@umich.edu}
\date{\today}
\thanks{The author was supported by NSF grant DMS-0502170.  This paper is part of the author's Ph.D. thesis at the University of Michigan under the direction of J.T. Stafford.}

\begin{document}
\begin{abstract}
Given $\zed$-graded rings $A$ and $B$, we ask when the graded module categories $\rgr A$ and $\rgr B$ are equivalent.  Using $\zed$-algebras, we relate the Morita-type results of 
\'Ahn-M\'arki and del R\'io 
to the twisting systems  introduced by Zhang, and prove, for example:

\medskip
\noindent{\bf Theorem.}
{\it If  $A$ and $B$ are $\zed$-graded rings, then:

(1)  $A$ is isomorphic to a Zhang twist of $B$ if and only if the $\zed$-algebras $\bbar{A} = \bigoplus_{i,j \in \zed} A_{j-i}$ and $\bbar{B} = \bigoplus_{i,j \in \zed} B_{j-i}$ are isomorphic.  

(2) If $A$ and $B$ are connected graded with $A_1 \neq 0$, then 
$\rgr A \simeq \rgr B$ if and only if $\bbar{A}$ and $ \bbar{B}$ are isomorphic.}
\medskip

This simplifies and extends Zhang's results.  
\end{abstract}
\maketitle

\section{Introduction}

The subject of this paper is graded Morita theory and its applications.  Given graded rings $A$ and $B$, algebras over some commutative ring $k$, we seek to find necessary and sufficient conditions for the categories $\rgr A$ and $\rgr B$ of graded right $A$ and $B$-modules to be equivalent; we write
 $\rgr A \simeq \rgr B$.   In the body of the paper we consider rings graded by arbitrary groups; for now we assume that $A$ and $B$ are $\zed$-graded.

In 1996, James Zhang \cite{Zh} solved the graded Morita problem for connected graded rings (i.e., 
$B = \bigoplus_{i \geq 0} B_i$ and $B_0 = k$, a field).  He defined  a {\em twisting system} $\tau = \{ \tau_n \}$ on a graded ring $B$, and used $\tau$ to define a new ``twisted'' multiplication on $B$.  This new ring is denoted $B^\tau$, and  is now known as a {\em Zhang twist} of $B$.  Zhang showed that if $\tau$ is a twisting system, then $\rgr  B \simeq \rgr  B^\tau$.  Further, he proved that if $A$ and $B$ are connected graded $k$-algebras with $A_1 \neq 0$, then $\rgr  A \simeq \rgr  B$ if and only if  $A \cong B^\tau$ for some twisting system $\tau$ on $B$.  

There are a number of general results establishing Morita-type theorems for graded module categories; see in particular the work of Angel del R\'io \cite{dR, dR2}.    Ultimately, del R\'io's results can be viewed as a special case of the  Morita theory for  rings with local units developed in \cite{Ab, AM}, which also  specializes to give the classical Morita theorems in the case where the grading group is trivial.  However, 
Zhang's results on twisting systems are formally quite different from classical Morita theory and from the results in \cite{dR, dR2, Ab, AM}.  The main aim of this paper is to unify the two theories. 

 Our technique involves {\em $\zed$-algebras}  --- roughly speaking, infinite matrix rings associated to graded rings.   More precisely,  a $\zed$-algebra $R = \bigoplus_{i,j\in \zed}R_{ij}$ is  a ring without 1 satisfying $R_{ij} R_{jl} \subseteq R_{il}$, with $R_{ij} R_{ml} = 0$ if $j \neq m$, such that each subring $R_{ii}$ contains a unit $1_{i}$.    In particular, given a  graded ring $A = \bigoplus_{i \in \zed} A_i$, then the ring $\bbar{A} = \bigoplus_{i,j \in \zed} A_{j-i}$ is a $\zed$-algebra with $\bbar{A}_{ij} = A_{j-i}$ and multiplication induced from $A$; we call it the {\em $\zed$-algebra associated to $A$}.  Precise definitions are in Section~\ref{sec-Galg}.

If $A$ and $B$ are $\zed$-graded rings, then by applying the Morita-type results mentioned above to the $\zed$-algebras $\bbar{A}$ and $\bbar{B}$, we can make the somewhat awkward definition of a Zhang twist (given in Section~\ref{sec-twist}) look very natural and   give an easy proof of Zhang's main theorem.   To do this, we recall that graded module categories have canonical automorphisms given by shifting degrees.  If $A$ is a $\zed$-graded ring and $M$ is  a graded right $A$-module, we denote the $n$'th shift of $M$ by $M\ang{n}$; that is, $M\ang{n} = \bigoplus_{i \in \zed} M\ang{n}_i$ is the graded right $A$-module given by $M\ang{n}_i = M_{i-n}$.  

We obtain:

\begin{theorem}\label{ithm0} 
{\em (Theorem~\ref{prop-TFAE}, Proposition~\ref{prop-cg})}
Let $A$ and $B$ be $\zed$-graded $k$-algebras.  Then:

(1) The associated $\zed$-algebras $\bbar{A}$ and $\bbar{B}$ are isomorphic if and only if there is an equivalence of categories $\Phi: \rgr A \to \rgr B$ such that $\Phi(A\ang{n}) \cong B \ang{n}$ for all $n \in \zed$.

(2) If $A$ and $B$ are connected graded and $A_1 \neq 0$, then $\rgr A \simeq \rgr B$ if and only if $\bbar{A}$ and $\bbar{B}$ are isomorphic.
\end{theorem}

\begin{theorem}\label{ithm1}
{\em (Corollary~\ref{cor-TFAE})}
Let $A$ and $B$ be  $\zed$-graded $k$-algebras.  Then   $B$ is a Zhang twist of $A$ if and only if  the associated $\zed$-algebras $\bbar{A}$ and $\bbar{B}$ are isomorphic.  

\end{theorem}

Thus we immediately obtain Zhang's main result:
\begin{corollary}\label{icor}
Let $A$ and $B$ be $\zed$-graded $k$-algebras.  Then:

(1) $B$ is isomorphic to a Zhang twist of $A$ if and only if there is an equivalence of categories $\Phi: \rgr A \to \rgr B$ such that $\Phi(A\ang{n}) \cong B \ang{n}$ for all $n \in \zed$.

(2) If $A$ and $B$ are connected graded and $A_1 \neq 0$, then $\rgr A \simeq \rgr B$ if and only if $B$ is isomorphic to a Zhang twist of $A$.
\end{corollary} 

Not all properties preserved under Morita equivalence are preserved under equivalences of graded categories --- notably, \cite{Zh} gives an example of $\zed$-graded rings $A$ and  $B$ with $\rgr A \simeq \rgr B$, where $A$ is gr-simple and prime and $B$ is neither.  However, in this example $B$ is semiprime.    It is not known whether, if $A$ and $B$ are $\zed$-graded rings such that  $\rgr A \simeq \rgr B$ and $A$ is semiprime, then $B$ must also  be semiprime, although \cite{Zh} shows it for connected graded Goldie rings.  As a step towards this general question, we prove:

\begin{proposition}\label{ithm-domain}
{\em (Proposition~\ref{thm-domain})}
Let $A$ and $B$ be $\zed$-graded $k$-algebras with $\rgr A \simeq \rgr B$.  If $B$ is a right Ore domain that is {\em fully graded} (that is, the grading on $A$ does not restrict to any proper subgroup of $\zed$), and $A_A$ has uniform dimension 1, then $A$ is a right Ore domain.  Further, if 
$Q$ is the graded quotient ring of $B$ and $Q'$ is the graded quotient ring of $A$, then $Q_0 \cong Q'_0$.
\end{proposition}

We also give an example showing that the property of having a graded semisimple graded Artinian graded quotient ring is not in general preserved under equivalences of graded module categories.

Using $\zed$-algebras to understand module categories has proven useful in other contexts; see, for example, \cite{BGS, SV}, among others.  We have  also used related techniques to study the category of graded modules over  the first Weyl algebra, $A_1$, under the Euler gradation \cite{S-A1}.   In that paper, we  completely classify graded rings $B$ such that $\rgr B$ is equivalent to $\rgr A_1$.   There are some surprising examples; in particular, there is a ring $S$ that is an idealizer in a localization of $A_1$ such that   $\rgr A_1 \simeq \rgr S$.   This is unexpected, not least because $A_1$ and its localizations are simple, whereas the idealizer has a finite dimensional representation and is not a maximal order.  

{\bf Notation.}  Throughout, we fix  a commutative ring $k$.  
If $G$ is a group and $A$ is a $G$-graded $k$-algebra, the category of ungraded right $A$-modules is denoted $\rmod A$, and the category of $G$-graded right $A$-modules is denoted $\rgr  A$.    If $M,N$ are objects of $\rgr  A$, then $\Hom_{\rgr  A}(M,N)$ is the set of graded homomorphisms of degree $e$; that is, the set of maps $\phi \in \Hom_{\rmod A}(M,N)$ so that $ \phi(M_g) \subseteq N_g$ for all $g \in G$.
We write:
\[ \hom_A(M,N)= \Hom_{\rgr A} (M,N)\]
and
\[ \Hom_A(M, N) = \Hom_{\rmod A} (M, N).\]
If $A$ is a ring without 1, we will denote the full subcategory of $\rmod A$ consisting of {\em unitary} modules --- that is, modules $M$ such that $M A = M$ --- by $\rmodu A$.
We similarly define $A \lgr$, $A \lmod$, and $A \lmodu$.

Given $g \in G$ and a right $A$-module $M = \bigoplus_{h \in G} M_h$ in $\rgr  A$, we define the  {\em $g$-th shift of $M$} to be a new object  $M \ang{g} = \bigoplus_{h \in G} M \ang{g}_h$ in $\rgr  A$, given by  $M\ang{g}_h = M_{g^{-1} h}$.   If $N = \bigoplus_{h \in G} N_h$ is a left $A$-module, then we define $N \ang{g}$ via 
$N \ang{g}_h = N_{h g^{-1}}$.
The {\em $g$-th shift functor} is the automorphism of the category $\rgr A$ (or $A \lgr$) sending $M$ to $M \ang{g}$.

\section{$G$-algebras and graded Morita theory for group-graded rings}\label{sec-Galg}

Before proceeding to the main results of this paper, we must formally introduce $G$-algebras. In this section, we give precise definitions, and show how to use the $G$-algebras associated to graded rings to apply the results of \cite{AM} to graded module categories.  This gives an extremely concrete solution of the Morita problem for graded rings, which we will use in the sequel.

Fix a group $G$.  Following \cite{BGS}, we define  
a {\em $G$-algebra} $R$ to be a $k$-algebra (possibly without 1) such that as a $k$-module $R = \bigoplus_{f, g \in G} R_{f,g}$, with  multiplication occuring matrix-wise:  that is, for all $f,g,h,g' \in G$, we have $R_{f,g} R_{g,h} \subseteq R_{f,h}$, and $R_{f,g} R_{g',h} = 0$ if $g \neq g'$.  We require that each ``diagonal'' subring $R_{f,f}$ have a unit $1_f$ that acts as a right identity on each $R_{g,f}$ and a left identity on each $R_{f,g}$.

We say that  $G$-algebras  $R$ and $S$ are {\em isomorphic as $G$-algebras} if there is a $k$-algebra  isomorphism $\phi: R \to S$ such that  $\phi(1_f) =  1_f$ for all $f \in G$.  More generally, we say that a $k$-algebra map (or $k$-module map) $\psi: R \to S$ is {\em graded of (left) degree $h$} if $\psi(R_{f,g}) \subseteq S_{hf,hg}$ for all $f,g \in G$.  

As \cite[Section~11]{SV} points out in the case $G=\zed$,  a  $G$-algebra can be viewed as a generalization of a graded ring in the following way:  if $A$ is a $G$-graded ring, define a $G$-algebra $\bbar{A} = \bigoplus_{f,g \in G} \bbar{A}_{f,g}$, the {\em (right) $G$-algebra associated to $A$}, by $\bbar{A}_{f,g} = A_{f^{-1}g}$.  Then  $\bbar{A}$ is a $G$-algebra under the  multiplication induced from $A$, since $\bbar{A}_{f,g} \cdot \bbar{A}_{g,h} = A_{f^{-1}g} \cdot A_{g^{-1}h} \subseteq A_{f^{-1} h} = \bbar{A}_{f,h}$.  There is a canonical isomorphism between the right module  $ A \ang{g}$  and the row $\bigoplus_{h \in G} \bbar{A}_{g,h}$  for any $g \in G$.  
Any $G$-graded right $A$-module is naturally a unitary right $\bbar{A}$-module; likewise any unitary right $\bbar{A}$-module  has a natural graded $A$-module structure given by $M_f = M \cdot 1_f$.  Thus we have an equivalence (in fact, an isomorphism) of categories between $\rgr A$ and $\rmodu \bbar{A}$.  In the sequel, we will identify the two categories.
(There is also an isomorphism between $A \lgr$ and $\bbar{A} \lmodu$, which we will not use.)

For completeness, we give the equivalent construction on the left:  define  the {\em left $G$-algebra associated to $A$} to be  $\hat{A} = \bigoplus_{f,g \in G} \hat{A}_{f,g}$, where $\hat{A}_{f,g} = A_{fg^{-1}}$.  As before, $A \lgr = \hat{A} \lmodu$, and the column $\bigoplus_{h \in G} \hat{A}_{h,g}$ is naturally isomorphic to the left module $A \ang{g}$.  For an alternate construction of $\hat{A}$, see the smash product constructed in \cite{B}.

There are $G$-algebra versions of the $\Hom$ functors on module categories.  If $R$ is a $G$-algebra and $S$ is a $G'$-algebra, then a unitary $(R, S)$-bimodule ${}_RQ_S$ has a natural bigraded structure $ Q = \bigoplus_{g \in G, h \in G'} Q_{g,h}$ 
given by $Q_{g,h} = 1_g \cdot Q \cdot 1_h$.  For such a bigraded module, we will denote the row $\bigoplus_{h \in G'} Q_{g,h}$ by $Q_{g*}$.

The bigraded bimodule ${}_RQ_S$ defines a covariant functor
$\Hfunct{S}{Q}{\blank}: \rmodu S \to \rmodu R$
given by 
\[ \Hfunct{S}{Q}{N} = \bigoplus_{g \in G} \Hom_S(Q_{g*}, N).\]
We leave to the reader the verification that $\Hfunct{S}{Q}{N}$ is a submodule of the right $R$-module $\Hom_S(Q, N) = \prod_{g \in G} \Hom_S(Q_{g*}, N)$, that $\Hfunct{S}{Q}{N} \cdot 1_g = \Hom_S(Q_{g*}, N)$, and that $\Hfunct{S}{Q}{N}$ is  the largest unitary $R$-submodule of $\Hom_S(Q,N)$; that is, $\Hfunct{S}{Q}{N} = \Hom_S(Q, N) \cdot R$.

In particular, if $B = \bigoplus_{h \in G'} B_h$ is a $G'$-graded ring and $M = \bigoplus_{g \in G, h \in G'} M_{g,h}$ is a $(G, G')$-bigraded right $B$-module (i.e., each $M_{g*}$ is a graded right $B$-module), then $M$ induces a functor $\Hfunct{\bbar{B}}{M}{\blank}$.  If $M$ is {\em locally finite} --- i.e., each $M_{g*}$ is a finitely generated $B$-module --- then we have $\Hfunct{\bbar{B}}{M}{M} = \bigoplus_{f, g \in G} \Hom_{\bbar{B}}(M_{g*}, M_{f*})$, and $\Hfunct{\bbar{B}}{M}{M}$ becomes a $G$-algebra, with 
\[\Hfunct{\bbar{B}}{M}{M}_{f,g} = \Hom_{\bbar{B}}(M_{g*}, M_{f*}) = \hom_B (M_{g*}, M_{f*});\] multiplication is given by composition of functions.  We refer to $\Hfunct{\bbar{B}}{M}{M}$ as the {\em endomorphism $G$-algebra of} $M$.  In this setting, $M$ is a bimodule over $\bbar{B}$ and its endomorphism $G$-algebra and so $\Hfunct{\bbar{B}}{M}{\blank}$ is a covariant functor from $\rgr B$ to $\rmodu \Hfunct{\bbar{B}}{M}{M}$.

If $R$ is a $G$-algebra, it has no unit unless $G$ is finite; however, $R$ is a {\em ring with local units} in the sense of \cite{Ab, AM}; that is, for any finite subset $X$ of $R$, there is an idempotent $a \in R$ such that $ax = x = xa$ for all $x \in X$.     Abrams \cite{Ab} and \'Ahn-M\'arki  \cite{AM} have generalized the classical Morita theorems to rings with local units.  In the next proposition, we reframe these results to analyze Morita theory for graded module categories in terms of $G$-algebras.   Similar results were obtained in \cite{dR, dR2}.  

We note that most results of \cite{AM} are stated for left module categories; however, by \cite[Corollary~2.3]{AM}, symmetric results hold on the right.

\begin{proposition}\label{thm-grMor}
Let $G$ and $G'$ be groups and let $A = \bigoplus_{g \in G} A_g$ and $B = \bigoplus_{h \in G'} B_{h}$ be graded rings.  Then $\rgr A \simeq \rgr B$ if and only if there is a $(G, G')$-bigraded right $B$-module $P_B = \bigoplus_{g \in G, h \in G'} P_{g, h}$ such that 

(1) $P_B$ is a locally finite projective generator for $\rgr B$; 

(2) $\bbar{A} \cong \Hfunct{\bbar{B}}{P}{P}$ as $G$-algebras.  

Further, if $P$ is as above, then $\Hfunct{\bbar{B}}{P}{\blank}: \rgr B \to \rgr A$ and $\blank \otimes_{\bbar{A}} P: \rgr A \to \rgr B$ are inverse equivalences; and if $\Phi: \rgr A \to \rgr B$ is an equivalence of categories, then $P = \Phi(\bbar{A}_A)$ satisfies (1) and (2), and $\Phi \cong \blank \otimes_{\bbar{A}} P$.
\end{proposition}
\begin{proof}
Suppose that $\Phi: \rgr A \to \rgr B$ is an equivalence of categories.  We define a $(G, G')$-bigraded right $B$-module $P = \bigoplus_{g \in G, h \in G'} P_{g,h} $, where $P_{g,h} = \Phi(A \ang{g})_h$.   Then $P_{g*} =\bigoplus_{h \in G'} P_{g,h} = \Phi(A \ang{g})$, and $P \cong \Phi(\bbar{A}_A)$.  Because $\bbar{A}_A$ is a locally finite projective generator, so is $P_B$; and functoriality of $\Phi$ gives $\Hfunct{\bbar{B}}{P}{P} = \Hfunct{\bbar{B}}{\Phi(\bbar{A})}{\Phi(\bbar{A})} \cong \Hfunct{\bbar{A}}{\bbar{A}}{\bbar{A}} = \bbar{A}$.  This direction also follows from the right-handed version of \cite[Theorem~2.1]{AM}, as does the fact that $\Phi \cong \blank \otimes_{\bbar{A}} P$.

Now suppose $P_B$ is a $(G, G')$-bigraded module satisfying (1) and (2).   Then  for any $N \in \rgr B$, $\Hfunct{\bbar{B}}{P}{N}$ is naturally isomorphic to $\Hom_{\bbar{B}}(P, N) \cdot \bbar{A}$.
 For any finite subset $F \subseteq G$, let $1_F = \sum_{f \in F} 1_f$ be the associated idempotent in $\bbar{A}$.  Clearly $ \varinjlim \Hom_{\bbar{B}}(1_F \cdot P, 1_F \cdot P) \cong \Hfunct{\bbar{B}}{P}{P} \cong \bbar{A}$.  Now we may apply the right-handed versions of \cite[Theorem~2.4, Theorem~2.5]{AM}, which together say that if $P_{\bbar{B}}$ is a locally finite projective generator  with $ \varinjlim \Hom_{\bbar{B}}(1_F \cdot P, 1_F \cdot P) \cong S$ for some ring with local units $S$, then $ \Hom_{\bbar{B}}(P, \blank) \cdot S: \rmodu \bbar{B} \to \rmodu S$ and $\blank \otimes_{S} P: \rmodu S \to \rmodu \bbar{B}$ are inverse equivalences.  Putting $S = \bbar{A}$, we have proved that $\Hfunct{\bbar{B}}{P}{\blank} = \Hom_{\bbar{B}}(P, \blank) \cdot \bbar{A}$ and $\blank \otimes_{\bbar{A}} P$ are inverse equivalences.  
\end{proof}

We remark that 
not all locally finite projective generators define an equivalence of graded module categories; that is, condition (2) above is nontrivial.  We give an example in  Example \ref{eg-notprincipal}.

\section{Principal $G$-algebras}\label{sec-principal}

We say that a $G$-algebra $R$ is {\em principal} if there is a $G$-graded ring $A$ such that $R \cong \bbar{A}$ as $G$-algebras.    (Of course, we then have $\rmodu R \simeq \rgr A$.)   
Since by Proposition~\ref{thm-grMor} equivalences of graded module categories involve principal $G$-algebras,   we are naturally interested in understanding these $G$-algebras better.    
In this section, we give a criterion for a $G$-algebra to be principal, and show that  $\bbar{A}$ and $ \bbar{B}$ are isomorphic as $G$-algebras exactly when $\rgr A$ and $\rgr B$ are related by an equivalence of a particularly nice form.

If  $R = \bbar{A}$ for some $G$-graded ring $A$, we make the straightforward observation that  for each $g \in G$ there is a map 
\begin{equation}\label{eq-shift}
 \shift_g: R = \bigoplus_{h,l \in G} R_{h,l} \to R = \bigoplus_{h, l \in G} R_{gh,gl}
 \end{equation}
defined on the component $R_{h, l}$ via the identifications
$R_{h,l} = A_{h^{-1}l} = A_{h^{-1}g^{-1} g l} =R_{gh,gl}$. 
Each $\shift_g$ is a graded $G$-algebra automorphism of degree $g$, and clearly $\shift_g \circ \shift_h = \shift_{gh}$.  The maps $\shift_g$ relate the  multiplications $\m_A$ in $A$ and $\m_R$ in $R$ as follows:  
if $x \in A_g = R_{e,g}$ and $y \in A_h = R_{e,h}$, then
\begin{equation}\label{eq-fund}
x \m_A y = x \m_R \shift_g(y) 
\end{equation}
in $R_{e,g} \m_R R_{g, gh} \subseteq R_{e,gh} = A_{gh}$.

The next proposition shows that in fact the existence of the maps $\shift_g$ characterizes principal $G$-algebras.    Given a $G$-algebra $R$, we will write $\Gaut(R)$ to mean the set of graded $k$-algebra  automorphisms of $R$ that have left degree $g$ for some $g \in G$.

\begin{proposition}\label{prop-principal}
If $R$ is a $G$-algebra, then $R$ is principal if and only if there is a group monomorphism
\begin{align*}
\T: G	& \hookrightarrow \Gaut(R) \\
	g & \mapsto \T_g
\end{align*}
such that for all $g$, $\T_g$ has left degree $g$.  

In particular, a $\zed$-algebra $R$ is principal if and only if $R$ has a $k$-algebra automorphism of degree $1$.
\end{proposition}

\begin{proof}
We have already seen that if $R$ is a principal $G$-algebra, then  there is such a map.

Conversely, suppose we are given $G$, $R$, and $\T$ as above.  We define a $G$-graded ring $A = \bigoplus_{g\in G}A_g$ as follows:   first define $A$ as a graded $k$-vector space via $A_g = R_{e,g}$.  Let $x \in A_g$, $y \in A_h$.  Mimicking \eqref{eq-fund}, we define the multiplication $\m_A$ on $A$ by
\begin{equation}\label{eq-star}
x \m_A y = x \m_R \T_g(y).
\end{equation}
Since $\T_g(y) \in R_{g,gh}$, we have 
$x \m_A y \in  R_{e, gh} = A_{gh}$ as required.

$A$ is easily seen to be $G$-graded, with unit $1 = 1_e$.  We check associativity of $\m_A$.    Let $z \in A_f$.   We have:
\[
 (x \m_A y ) \m_A z = (x \m_R \T_g(y)) \m_R \T_{gh}(z)  =  x \m_R \T_g(y) \m_R \T_g \T_h(z) 
\]
since $\T$ is a group homomorphism.
But each $\T_g$ is a $G$-algebra homomorphism, so 
\[ x \m_A (y \m_A z) = x \m_R \T_g(y \m_R \T_h(z)) = x \m_R \T_g(y) \m_R \T_g \T_h(z)\]
and we have $(x \m_A y) \m_A z   = x \m_A (y \m_A z)$.    Thus $A$ is a $G$-graded $k$-algebra.

We show that $R \cong {\bbar{A}}$ as $G$-algebras.  Define $\alpha: {\bbar{A}} \to R$ by letting
$\alpha$ act on ${\bbar{A}}_{f,g}$ via:
\[ \xymatrix{
{\bbar{A}}_{f,g} \ar[r]^(.25){\shift_f^{-1}} & {\bbar{A}}_{e,f^{-1}g} = A_{f^{-1}g} = R_{e,f^{-1}g} \ar[r]^(.75){\T_f}	& R_{f,g}.	} \]

This is a $k$-linear bijection, and if $x \in {\bbar{A}}_{f,g}$, $y \in {\bbar{A}}_{g,h}$,  using  \eqref{eq-fund} and \eqref{eq-star}, we have
\begin{multline*}
\alpha(x \m_{\bbar{A}} y) = \T_f \shift_f^{-1} (x \m_{\bbar{A}} y) = \T_f \bigr(\shift_f^{-1}(x) \m_{\bbar{A}} \shift_f^{-1}(y) \bigl) \\
	= \T_f\bigl(\shift_f^{-1} (x) \m_A \shift^{-1}_g(y)\bigr)
	= \T_f\bigl(\shift^{-1}_f (x) \m_R \T_{f^{-1}g} \shift^{-1}_g(y)\bigr) \\ = \T_f \shift_f^{-1} (x) \m_R \T_g \shift_g^{-1}(y) = \alpha(x) \m_R \alpha(y).
\end{multline*}
Thus $R \cong \bbar{A}$ is principal.
\end{proof}

If $R, G, \T$ are as in Proposition \ref{prop-principal}, we say that $\T$ is a {\em principal map of $R$}.   We call the ring $A$  the {\em compression of $R$ by $\T$} and we  write $A = R^{\T}$.   If $R = \bbar{A}$, we say that the principal map $\shift$ defined in \eqref{eq-shift} is the {\em canonical principal map of $\bbar{A}$}. 

\begin{example}\label{eg-notprincipal}
We give an example of a locally finite projective generator whose endomorphism $G$-algebra is not principal, showing that condition (2) of Proposition~\ref{thm-grMor} is nontrivial.  Let $G = G' = \zed$, $B = k[x]$ for some field $k$, and let
\[
P_{n*} = \left\{ \begin{array}{ll}
		B\ang{n}	& \mbox{if $n$ is even} \\
		B\ang{-n}	& \mbox{if $n$ is odd.}  \end{array} \right. \]
Then $P$ is a locally finite projective generator for $\rgr B$.  We   show that $H = \Hfunct{\bbar{B}}{P}{P}$ is not a principal $G$-algebra. We have: 
\[H_{0,2} = \hom_B(B\ang{2}, B) = B_2 = k x^2,\]
 but 
\[H_{1,3} = \hom_B(B\ang{-3}, B\ang{-1}) = B_{-2} = 0. \]
 Thus $H$ has no principal map, and by   Proposition \ref{prop-principal}, $H$ is not isomorphic to $\bbar{A}$ for any $\zed$-graded ring $A$.  
\end{example}

Let $R \cong \bbar{A}$ be a principal $G$-algebra.  If we also have $R \cong \bbar{B}$ for some other $G$-graded ring $B$, then obviously the categories $\rgr A = \rgr \bbar{A}$ and $\rgr B = \rgr \bbar{B}$ are isomorphic.  In fact, as the next result shows, this isomorphism is an equivalence of a particularly nice form; furthermore, the existence of such a ``nice'' equivalence between $\rgr A$ and $\rgr B$ implies that the associated $G$-algebras $\bbar{A}$ and $ \bbar{B}$ are isomorphic.  This proves part (1) of Theorem~\ref{ithm0}.

\begin{theorem}\label{prop-TFAE}
Let $A$ and $B$ be $G$-graded $k$-algebras.  The following are equivalent:

(1)  For some principal map $\T$ of $\bbar{B}$, $\bbar{B}^{\T} \cong A$ via a degree-preserving isomorphism.

(2) The right associated $G$-algebras $\bbar{A}$ and $\bbar{B}$ are isomorphic as $G$-algebras.

(3) There is an equivalence of categories  $\Phi: \rgr  A \to \rgr  B$ such that $\Phi(A\ang{g}) \cong B \ang{g}$ for all $g \in G$.

(2') The left associated $G$-algebras $\hat{A}$ and $\hat{B}$ are isomorphic as $G$-algebras.

(3') There is an equivalence of categories  $\Phi':  A \lgr \to  B \lgr$ such that $\Phi'(A\ang{g}) \cong B \ang{g}$ for all $g \in G$.
\end{theorem}

\begin{proof}
(3) $\Rightarrow$ (2).  Since $\Phi(A\ang{g}) \cong B \ang{g}$ for all $g$, we have $\Phi(\bbar{A}_A) \cong \bbar{B}_B$.  Because $\Phi$ is a category equivalence, Proposition~\ref{thm-grMor} gives us
an isomorphism of $G$-algebras between $\bbar{A}$ and $\Hfunct{\bbar{B}}{\Phi(\bbar{A})}{\Phi(\bbar{A})}$.  This is isomorphic to $\Hfunct{\bbar{B}}{\bbar{B}}{\bbar{B}} = \bbar{B}$.

(2) $\Rightarrow$ (3).  As the identification between $\rgr  B$ and $\rmodu  \bbar{B}$ commutes with shifting,  it is enough to show that there is an equivalence $\Phi: \rmodu  \bbar{A} \to \rmodu \bbar{B}$ such that $\Phi(A\ang{g}) \cong B \ang{g}$ for all $g \in G$.  But this is clear, as the isomorphism between $\bbar{A}$ and $\bbar{B}$ preserves degree.

(1) $\Rightarrow$ (2).   The proof of Proposition \ref{prop-principal} shows that $\bbar{B} \cong \bbar{(\bbar{B}^{\T})} \cong \bbar{A}$.

(2) $\Rightarrow$ (1).  Let $\alpha: \bbar{A} \to \bbar{B}$ be an isomorphism of $G$-algebras, and let $\shift$ be the canonical principal map of $\bbar{A}$.   For all $g \in G$, define $\T_g = \alpha \shift_g \alpha^{-1}$.  Then clearly $\T$ is a principal map of $\bbar{B}$, and $\alpha$ gives a graded $k$-linear bijection from $A = \bbar{A}^{\shift} \to \bbar{B}^{\T}$.  We check this is a ring homomorphism: since $\alpha$ preserves degree, if $a \in A_g$ and $b \in A_h$, then we have
\begin{align*}
\alpha(a \m_A b) = \alpha(a \m_{\bbar{A}} \shift_g(b)) & = \alpha(a) \m_{\bbar{B}} \alpha( \shift_g (b)) = \alpha(a) \m_{\bbar{B}} \T_g (\alpha(b)) = \alpha(a) \bullet \alpha(b) 
	\end{align*}
where $\bullet$ denotes multiplication in $\bbar{B}^{\T}$.

$(2')$ $\iff$ $(3')$ follows by symmetry, and (2) $\iff$ $(2')$ is obvious from the definitions of $\bbar{A}$ and $\hat{A}$.
\end{proof}

{\bf Remark:} There are results describing the conditions under which equivalences of graded module categories commute with shifting; see \cite{BR, GG}.  We note that condition (2) of Theorem~\ref{prop-TFAE} is much weaker than the condition that $\Phi$ and shifting by $g$ commute as functors; in fact (2) does not even imply that $\Phi(M \ang{g}) \cong (\Phi(M))\ang{g}$ for all $M$.

To complete this section, we specialize to connected graded rings and show that here all equivalences are of the form given in Theorem~\ref{prop-TFAE}.  This gives part (2) of Theorem~\ref{ithm0}.   Recall that if $k$ is a field, a $\zed$-graded $k$-algebra $A$ is {\em connected graded} if $A = \bigoplus_{i \geq 0} A_i$ and $A_0 = k$.  

\begin{proposition}\label{prop-cg}
Let $A$ and $B$ be $\zed$-graded and connected graded $k$-algebras with $A_1 \neq 0$.  Then $\rgr A \simeq \rgr B$ if and only if $A$ and $B$ satisfy the equivalent conditions of Theorem~\ref{prop-TFAE}.
\end{proposition}
\begin{proof}

We follow the approach of \cite[Theorem~3.5]{Zh}.  Suppose  
$\Phi: \rgr A \to \rgr B$ is an equivalence.  By Proposition~\ref{thm-grMor},  without loss of generality we are in the following situation:   $P_B = \Phi(\bbar{A})$ is a $(\zed, \zed)$-bigraded projective generator for $\rgr B$; $\Hfunct{\bbar{B}}{P}{P} \cong \bbar{A}$; and $\Phi = \blank \otimes_{\bbar{A}} P$.  

Each $P_{n*} = \Phi(A\ang{n})$ is an indecomposable projective in $\rgr B$.  As $B$ is connected graded, there is some integer $f(n)$ such that $P_{n*} = B \ang{f(n)}$.  Since for all $m \in \zed$, $\Hfunct{\bbar{B}}{P\ang{m}}{P \ang{m}} \cong \Hfunct{\bbar{B}}{P}{P}$, we may shift $P$  so that without loss of generality $f(0) = 0$.

Since $P$ generates $\rgr B$, $f$ is surjective.  Because $A$ and $B$ are both non-negatively graded, we have $f(m) = f(n) \iff P_{m*} \cong P_{n*} \iff A\ang{m}  \cong A \ang{n} \iff m = n$.  Thus $f: \zed \to \zed$ is bijective.  For all $n \in \zed$, 
\[B_{f(n+1)-f(n)} = \hom_B(B \ang{f(n+1)}, B \ang{f(n)}) \cong \Hfunct{\bbar{B}}{P}{P}_{n, n+1} \cong  \bbar{A}_{n, n+1}  = A_1 \neq 0.\]
Thus $f(n+1) - f(n) \geq 0$.  Since $f: \zed \to \zed$ is an increasing bijection with $f(0) = 0$, for all $n$, $f(n) = n$, and $\Phi(A \ang{n} ) \cong B\ang{n}$.  Thus the equivalent conditions of Theorem~\ref{prop-TFAE} are satisfied.  The other direction follows immediately from  Theorem~\ref{prop-TFAE}.  
\end{proof}

\section{$G$-algebras and twisting systems}\label{sec-twist}

In this section, we relate the results in \cite{Zh} to the material developed in Section \ref{sec-principal} and prove Theorem~\ref{ithm1}, thus unifying Zhang's results and classical Morita theory.  In particular, we will see that twisting systems naturally come from principal maps.

Throughout this section we fix a group $G$ and let $A = \bigoplus_{g \in G} A_g$ and $B = \bigoplus_{g\in G}B_g$ be $G$-graded $k$-algebras.  Recall \cite{Zh} that  a {\em twisting system} $\tau$ on $A$ is a set $\{\tau_g \left\lvert g \in G\right.\}$ of $k$-linear graded automorphisms of $A$ such that
\begin{equation}\label{eq-twist1}
\tau_g(y \tau_h (z)) = \tau_g(y) \tau_{gh}(z)
\end{equation}
for all $g,h,l \in G$, $y \in A_h$, and $z \in A_l$.  Using a twisting system $\tau$, the {\em twisted algebra} (now called a {\em Zhang twist}) $A^\tau$ is the graded $k$-module $A$, with new multiplication $\star$ defined by 
\[ x \star y = x \tau_g(y) \]
where $x \in A_g$.

The twisting relation \eqref{eq-twist1} is somewhat cumbersome. The next proposition shows that instead of twisting systems, we may equivalently study principal maps; in particular, if $G = \zed$, then a twisting system on $A$ is a set of maps $\{ \tau_n \}$ that come from powers of a single algebra automorphism of the $\zed$-algebra $\bbar{A}$.

\begin{proposition}\label{prop-twist-compress}
Let $A$ be a $G$-graded $k$-algebra, and let $R = \bbar{A}$ be the associated $G$-algebra.  There is a bijection between the set of principal maps $\T$ of $R$ and the set of 
twisting systems $\tau$ on $A$ in such a way that if $\tau$ and $\T$ correspond, then $A^\tau \cong R^{\T}$.

In particular, if $G = \zed$ then there is a bijection between twisting systems on $A$ and automorphisms of $\bbar{A}$ of degree 1.
\end{proposition}

\begin{proof}
Recall that $R$ has a canonical principal map $\shift$ defined by \eqref{eq-shift}.  
Define a function  $\Delta$ from the set of principal maps on $R$ to the set of twisting systems on $A$  as follows:  if $\T: G \to \Gaut(R)$ is a principal map, define $\tau = \Delta(\T)$ by 
\[ \tau_g = \shift_g^{-1} \T_g: R_{e*} = A \to A.\]
We check that $\tau$ is a twisting system.  Let $y \in A_h$, $z \in A_l$.  Then
\[
\tau_g(y \m_A \tau_h(z))	 = \shift_g^{-1} \T_g(y \m_A \shift_h^{-1} \T_h(z)) 
	 = \shift_g^{-1} \T_g(y \m_R \T_h(z)) 
\]
where the last equality is by \eqref{eq-fund}.
Also, we have
\begin{multline*}
\tau_g(y) \m_A \tau_{gh}(z) = \shift_g^{-1} \T_g(y) \m_A \shift_{gh}^{-1} \T_{gh} (z)  = \shift_g^{-1}\T_g(y) \m_R \shift_h \shift_{gh}^{-1}  \T_{gh}(z) \\
	 = \shift_g^{-1} \T_g(y) \m_R \shift_g^{-1} \T_g \T_h (z)  = \shift_g^{-1} \T_g(y \m_R \T_h(z))
	\end{multline*}
as $\shift_g$, $\T_g$ are $G$-algebra automorphisms.  Thus \eqref{eq-twist1} is satisfied, $\tau$ is a twisting system, and $\Delta$ is well-defined.

Conversely, define a map $\Gamma$ from the set of twisting systems on $A$ to the set of group homomorphisms from  $G$ to  $\Gaut(R)$ 
 as follows:  if $\tau$ is a twisting system, define $\Gamma(\tau)_g: R = \bigoplus_{h,l \in G} R_{h,l} \to R = \bigoplus_{h,l \in G} R_{gh, gl}$ as the map that on $R_{h,l}$ acts by 
\[\xymatrix{
 \Gamma(\tau)_g: R_{h,l} \ar[r]^(0.45){\shift_h^{-1}} & 
 	R_{e,h^{-1}l} = A_{h^{-1}l} \ar[r]^{\tau_{gh} \circ \tau_h^{-1}} &
	A_{h^{-1}l} = R_{e,h^{-1}l} \ar[r]^(0.6){\shift_{gh}} &
	R_{gh, gl}.		} \]
It is immediate from the definition that $\Gamma(\tau)_g$ is a $k$-linear bijection of degree $g$, and that $\Gamma(\tau)_g \circ \Gamma(\tau)_f = \Gamma(\tau)_{gf}$.  We postpone for the moment verifying that $\Gamma(\tau)_g$ is a ring homomorphism.

Assuming this, we claim that in fact $\Gamma = \Delta^{-1}$.   First, 
since on $R_{e*}$ we have $\Gamma(\tau)_g = \shift_g \circ \tau_g$, therefore
\[(\Delta\Gamma(\tau))_g =  \shift_g^{-1} \shift_g \tau_g = \tau_g,\]
 so $\Delta \Gamma = 1$.   Let  $\T$ be a principal map on $R$, and put $\tau = \Delta(\T)$.  Then $\Gamma\Delta(\T)_g$ acts on $R_{h,l}$ as:
\begin{align*}
 \shift_{gh} \tau_{gh} \tau_h^{-1} \shift_h^{-1} 	& = \shift_{gh} \shift_{gh}^{-1} \T_{gh} \T_h^{-1} \shift_h \shift_h^{-1} & \mbox{ by definition of $\Delta$} \\
	& = \T_g &  \mbox{ since $\T_{gh} = \T_g \T_h$.}
\end{align*}
Thus $\Gamma\Delta = 1$, $\Gamma$ and $\Delta$ are inverses, and in particular $\Delta$ is a bijection as claimed.

Fix a twisting system $\tau$ and let $\T = \Delta(\tau)$ be the corresponding principal map.  Let $\star$ denote multiplication in $A^\tau$ and $\bullet$ denote multiplication in $R^{\T}$.  We check that $\star$ and $\bullet$ are equal.  If $x \in (R^{\T})_g = R_{e,g} = A_g$ and $y \in (R^{\T})_h = A_h$, then
\[ x \bullet y = x \m_R \T_g(y)\]
by definition, and 
\begin{align*}
x \star y = x \m_A \tau_g(y) = x \m_A \shift_g^{-1} \T_g(y) & =  x \m_R \T_g(y) & \mbox{by \eqref{eq-fund}} \\
 & = x \bullet y .
 \end{align*}
  Thus $A^\tau \cong R^{\T}$.

To complete the proof, we must verify that if $\tau$ is a twisting system on $A$ and $\T = \Gamma(\tau)$, then each $\T_g$ is a ring homomorphism.  By \cite[Equation~2.1.3]{Zh}, for all $h, m, l \in G$, $a \in A_m$, and $b \in A_l$, we have
\[ \tau_h^{-1} (a \m_A b) = \tau_h^{-1}(a) \m_A \tau_m \tau_{hm}^{-1}(b).\]
Thus, applying \eqref{eq-twist1}, we have
\begin{equation}\label{eq-twist3}
\tau_{gh}(\tau_h^{-1}(a \m_A b)) = \tau_{gh}(\tau_h^{-1}(a)) \m_A  \tau_{ghm} (\tau_{hm}^{-1}(b)).
\end{equation}

Let $h,f,l \in G$, $x \in R_{h,f}$, and $y \in R_{f,l}$.  We have:
\begin{multline*} \T_g(x \m_R y) = \shift_{gh} \tau_{gh} \tau_h^{-1} \shift_h^{-1} (x \m_R y) 
	= \shift_{gh} \tau_{gh} \tau_h^{-1} \bigl(\shift_h^{-1} (x) \m_R \shift_h^{-1} (y) \bigr) \\
	= \shift_{gh} \tau_{gh} \tau_h^{-1} \bigl(\shift_h^{-1} (x) \m_A \shift_f^{-1} (y) \bigr) 
	\end{multline*}
by \eqref{eq-fund}.
Since $\shift_h^{-1}(x) \in R_{e,h^{-1}f} = A_{h^{-1}f}$, applying \eqref{eq-twist3} with $m = h^{-1}f$, we see that this is equal to 
\begin{multline*}
\shift_{gh}\bigl( \tau_{gh} \tau_h^{-1} \shift_h^{-1}(x) \m_A \tau_{gf} \tau_f^{-1} \shift_f^{-1}(y) \bigr) \\
\begin{array}{ll}
= \shift_{gh} \bigl( \tau_{gh} \tau_h^{-1}\shift_h^{-1} (x) \m_R \shift_{h^{-1}f} \tau_{gf} \tau_f^{-1} \shift_f^{-1}(y) \bigr)	& \mbox{by \eqref{eq-fund}} \\
= \bigl( \shift_{gh} \tau_{gh} \tau_h^{-1}\shift_h^{-1} (x) \bigr) \m_R \bigl( \shift_{gf}   \tau_{gf} \tau_f^{-1} \shift_f^{-1}(y) \bigr) & \\
= \T_g(x) \m_R \T_g(y) \\
\end{array} \end{multline*}
as required. \end{proof}

We thus may add another equivalent condition to Theorem~\ref{prop-TFAE}, proving part (1) of Theorem~\ref{ithm1}.

\begin{corollary}\label{cor-TFAE}
Let $A$ and $B$ be $G$-graded rings.  Then $A$ and $B$ satisfy the equivalent conditions of Theorem~\ref{prop-TFAE} if and only if $A$ is isomorphic to a twisted algebra of $B$.
\end{corollary}
\begin{proof}
Combine Proposition~\ref{prop-twist-compress} and Theorem~\ref{prop-TFAE}. \end{proof}

One part of Corollary~\ref{cor-TFAE} was proven by Zhang:
\begin{corollary}\label{cor-z-TFAE}
{\em (\cite[Theorem~3.3]{Zh})}
Let $G$ be  a group and let $A$ and $B$ be $G$-graded rings.  Then $A$ is isomorphic to a twisted algebra of $B$ if and only if there is an equivalence $\Phi: \rgr A \to \rgr B$ such that $\Phi(A\ang{g}) \cong B\ang{g}$ for all $g \in G$. \qed
\end{corollary}

Thus we also obtain Zhang's other  main result:
\begin{corollary}\label{z-cg}
{ \em (\cite[Theorem~3.5]{Zh})}
Let $k$ be a field, and let  $A$ and $B$ be two connnected graded and  $\zed$-graded $k$-algebras with $A_1 \neq 0$.  Then $B$ is isomorphic to a twisted algebra of $A$ if and only if $\rgr  A$ is equivalent to $\rgr  B$.
\end{corollary}
\begin{proof}
Combine Corollary \ref{cor-TFAE} and Proposition \ref{prop-cg}.
\end{proof}

\section{Graded domains}

In this section we apply Proposition~\ref{thm-grMor} to understand  {\em graded domains} --- i.e.  graded rings where all nonzero homogeneous elements are nonzerodivisors.  
For ungraded rings, being a domain is not a Morita invariant, although being prime is.  However, if $A$ is an ungraded Ore domain, and $B \cong \End_A(M, M)$ is Morita equivalent to $A$ and ``clearly not a matrix ring'' --- i.e. $B_B$ has uniform rank $1$ --- then  $B$ is also a right Ore domain and $Q(A) \cong Q(B)$.  Proposition \ref{thm-domain} is the graded analogue of this result.  

We say that a $G$-graded ring $B$ is {\em fully graded} if the set $\{g \in G \st B_g \neq 0\}$ is not contained in any proper subgroup of $G$.  

\begin{proposition}\label{thm-domain}
Let $A$ and $B$ be $G$-graded $k$-algebras with $\rgr A \simeq \rgr B$.  If $B$ is a fully graded right Ore domain and $A_A$ has uniform dimension 1, then $A$ is a graded right Ore domain.  Further, if $Q$ is the graded quotient ring of  $B$, and $Q'$ is the graded quotient ring of $A$, then $Q_e \cong Q'_e$.
\end{proposition}

If the grading group $G$ is ordered, then a graded domain is also an ungraded domain, and so this result implies Proposition~\ref{ithm-domain} from the Introduction.

\begin{proof}
Let $A$ and $B$ satisfy the hypotheses of the theorem, and let $Q$ be the graded quotient ring of $B$.  Then $Q$ is a graded division ring, so $Q_e = D$ is a division ring.  We note that $Q$ is the unique uniform graded injective torsion-free module in $\rgr B$, and so $Q \ang{g} \cong Q$ for all $g \in G$.   (The hypothesis that $B$ is fully graded is necessary to ensure this:  if $B = Q = k[x^2, x^{-2}]$ then the previous sentence is false.)   Further, $\bbar{Q}$ is a $(G,G)$-bigraded right $B$-module and a $G$-algebra, and the identifications $Q_{f^{-1}g} = \hom_B(Q\ang{g}, Q\ang{f})$ induce a canonical isomorphism of $G$-algebras between $\bbar{Q}$ and $\Hfunct{\bbar{B}}{\bbar{Q}}{\bbar{Q}}$.   In the $G$-algebra $\bbar{Q}$, each nonzero homogeneous element $q \in \bbar{Q}_{f,g}$ is invertible, in the sense that there is an element  $q^{-1} \in \bbar{Q}_{g,f}$, with $q\cdot q^{-1} = 1_f$ and  $q^{-1} \cdot q = 1_g$.

By Proposition~\ref{thm-grMor}, we may assume that we are in the following situation:  there is a category equivalence $\Phi: \rgr A \to \rgr B$, the $B$-module $P = \Phi(\bbar{A})$ is a locally finite $(G,G)$-bigraded projective generator for $\rgr B$, $\Phi = \blank \otimes_{\bbar{A}} P$, and $\Psi = \Hfunct{\bbar{B}}{P}{\blank}$ is the  inverse equivalence to $\Phi$.  Each $P_{g*} = \Phi(A\ang{g})$ is projective and therefore torsion-free, and uniform by hypothesis.  Thus for all $g \in G$, there is an injection $\gamma_g: P_{g*} \to Q \ang{g}$.  

Using the $\gamma_g$, we define a map $\Theta:  \Hfunct{\bbar{B}}{P}{P} \to \bbar{Q}$ as follows: let $f,g \in G$, $\phi  \in  \Hfunct{\bbar{B}}{P}{P}_{f,g} = \hom_B(P_{g*}, P_{f*})$.  As $Q \cong Q\ang{f}$ is graded injective, the map $\gamma_f \circ \phi: P_{g*} \to Q\ang{f}$ lifts uniquely via the inclusion $\gamma_g$ to a map $\bbar{\phi}: Q\ang{g} \to Q \ang{f}$ such that the diagram
\[ \xymatrix{
P_{g*} \ar[r]^{\phi} \ar[d]_{\gamma_g}	& P_{f*} \ar[d]^{\gamma_f} \\
Q\ang{g} \ar[r]_{\bbar{\phi}}	& Q \ang{f}  	} \]
commutes.  Put $\Theta(\phi) = \bbar{\phi} \in \Hfunct{\bbar{B}}{\bbar{Q}}{\bbar{Q}}$.

The map $\Theta$ is easily seen to be  a $G$-algebra injection from $\Hfunct{\bbar{B}}{P}{P} \cong \bbar{A}$ to  $\bbar{Q}$.   
To simplify notation, put $R = \Theta(\Hfunct{\bbar{B}}{P}{P}) \subseteq \bbar{Q}$.  Let $\T$ be the principal map of $R$ induced from the canonical principal map of $\bbar{A}$.   Thus $A \cong R^{\T}$.

The restriction of $\Theta$ to $\Hfunct{\bbar{B}}{P}{P}_{f*}$ gives an inclusion   $R_{f*} = \Theta(\Hfunct{\bbar{B}}{P}{P_{f*}}) \subseteq \bbar{Q}_{f*}$.  
Because $(\bbar{Q}_{f*})_R \cong \Hfunct{\bbar{B}}{P}{Q\ang{f}}$, it is a graded uniform right $R$-module.  Thus if $q \neq 0 \in \bbar{Q}_{f,g}$, then $qR \cap R_{f*} \neq 0$, so there are $h \in G$,  $r \in R_{f, h}$, $s \neq  0 \in R_{g,h}$ such that $qs = r$.  Thus $R$ is a graded right order in $\bbar{Q}$.   

For each $l \in G$, we use this to extend $\T_l: R \to R$ to an automorphism of $\bbar{Q}$ of degree $l$ by putting $\T_l(q) = \T_l(r) \T_l(s)^{-1} \in \bbar{Q}_{lf, lh} \cdot \bbar{Q}_{lh, lg} = \bbar{Q}_{lf, lg}$.  It is an easy verification that each $\T_l$ is well-defined and is a $k$-algebra automorphism of $\bbar{Q}$.  Compressing by $\T$, we have $A \cong R^{\T} \subseteq \bbar{Q}^{\T}$.  Put $Q' = \bbar{Q}^{\T}$.  

Now $Q'$ is a strongly $G$-graded graded division ring, and as $R$ is a right order in $\bbar{Q}$, $R^{\T} \cong A$ is clearly a right order in $Q'$.  Thus $A$ is a graded right Ore domain with graded quotient $Q'$.   Again, $D' = Q'_e$ is a division ring, and we have, by \cite[Theorem~A.I.3.4]{NV} and Proposition \ref{prop-principal}, 
\[ \rmod D' \simeq \rgr Q' \simeq \rmodu \bbar{Q} \simeq \rgr Q \simeq \rmod D.\]
Thus $D$ and $D'$ are Morita equivalent division rings, and therefore $D \cong D'$.
\end{proof}

In the ungraded situation, being prime, semiprime, and semiprime (right, left) Goldie are Morita invariants.  On the other hand, the following example from \cite{Zh} shows that being graded prime is not invariant under taking equivalences of graded module categories.  Let $R= k[x^2, x^{-2}]$ for some field $k$.  Then the rings
\[ A = \left( \begin{array}{ll}
		R	& x R \\
		xR	& R 
		\end{array} \right)
		\]
and
\[ B = k[x, x^{-1}] \oplus k [x, x^{-1}] \]
are Zhang twists of each other, so $\rgr A \simeq \rgr B$.  Here $A$ is a prime graded simple graded Artinian ring, and $B$ is graded semisimple graded Artinian and semiprime but not prime or even graded prime.

In general, we do not know whether being  graded semiprime graded Goldie (or graded semiprime) is invariant under equivalences of graded module categories.  \cite[Proposition~5.8]{Zh} shows that this is true for $\zed$-graded and connected graded rings.  The proof uses the fact that if $B$ is connected graded and $\zed$-graded, then  $B$ is semiprime graded right Goldie if and only if $B$ has a  graded semisimple graded Artinian right quotient ring, if and only if every graded essential right ideal of $B$ contains a homogeneous regular element.  This last property is invariant under Zhang twisting.  However, Example \ref{eg2} shows that in general, the property of being semiprime graded Goldie and having a graded semisimple graded Artinian quotient ring is not invariant under equivalences of graded module categories, even for $\zed$-graded rings.  Thus techniques that do not involve localizing must  be used  to answer the general question.

\begin{example}\label{eg2}
Let $k$ be a field, $G = \zed$, and let $A = k[x] \oplus k[y]$, with $x, y \in A_1$.  $A$ is semiprime graded Goldie and has a graded semisimple graded Artinian quotient ring. Define a bigraded module $P_A = \bigoplus_{n, m \in \zed} P_{nm}$ by $P_{n*} = k[x]\ang{n} \oplus k[y] \ang{-n}$.  Then $P$ is a locally finite projective generator for $\rgr A$, and if we put $H = \Hfunct{\bbar{A}}{P}{P}$ then we have
\[ H_{n,m} = \left\{ {\begin{array}{ll}
	k x^{n-m}	& \mbox{if $n > m$} \\
	k \oplus k	& \mbox{if $n = m$} \\
	k y^{m-n}	& \mbox{if $n< m$.}
	\end{array} } \right. \]
$H$ is clearly principal, and in fact $H \cong \bbar{B}$, where $B = k[x] \oplus k[y]$, with $x \in B_1$, $y \in B_{-1}$; thus $\rgr A \cong \rgr B$.  $B$ is also semiprime graded Goldie, but all homogeneous regular elements of $B$ are in degree 0, and are already invertible in $B$.  Thus $B$ does not have a graded semisimple graded Artinian quotient ring.
\end{example}

\section*{Acknowledgements}
We thank Michel van den Bergh for suggesting applying $\zed$-algebras to obtain the results in \cite{Zh}.

\end{document}